\documentclass[12pt, a4j]{article}
\AtBeginDvi{}
\usepackage{ascmac}
\usepackage{amsmath, amssymb, amscd}
\usepackage{type1cm}
\usepackage{bm}
\usepackage[pdftex,hiresbb]{graphicx}
\usepackage{float}
\usepackage{tikz}
\newtheorem{defi}{Definition}[section]
\newtheorem{rema}[defi]{Remark}
\newtheorem{exa}[defi]{Example}

\newtheorem{thm}[defi]{Theorem}

\makeatletter
    
    \@addtoreset{equation}{section}
  \makeatother

\begin{document}
\title{On equivariant index of a generalized Bott manifold}
\author{Yuki Sugiyama}
\date{}
\maketitle

{\small \textbf{Abstract.} \,  In this paper, we consider the equivariant index
of a generalized Bott manifold. We show the multiplicity function of the equivariant
index is given by
the density function of a generalized twisted cube.
In addition, we give a Demazure-type character formula of this representation.}

\section{\large Introduction}

A {\it Bott tower} of height $n$ is a sequence:
\begin{equation*}
M_n \stackrel{\pi_n}{\to} M_{n-1} \stackrel{\pi_{n-1}}{\to} \cdots
\stackrel{\pi_2}{\to} M_1 \stackrel{\pi_1}{\to} M_0 = \{ \textup{a point} \}
\end{equation*}
of complex manifolds $M_j = \mathbb{P}(\underline{\mathbb{C}} \oplus E_j)$,
where $\underline{\mathbb{C}}$ is the trivial line bundle over $M_{j-1}$,
$E_j$ is a holomorphic line bundle over $M_{j-1}$, $\mathbb{P}(\cdot)$
denotes the projectivization, and $\pi_j : M_j \to M_{j-1}$ is the projection of the
$\mathbb{C}P^1$-bundle. We call $M_j$ a $j$-{\it stage Bott manifold}.
The notion of a Bott tower was introduced by Grossberg and Karshon (\cite{GK}).

A {\it generalized Bott tower} is a generalization of a Bott tower.
A generalized Bott tower of height $m$ is a sequence:
\begin{equation*}
B_m \stackrel{\pi_m}{\to} B_{m-1} \stackrel{\pi_{m-1}}{\to} \cdots
\stackrel{\pi_2}{\to} B_1 \stackrel{\pi_1}{\to} B_0 = \{ \textup{a point} \},
\end{equation*}
of complex manifolds $B_j = \mathbb{P}(\underline{\mathbb{C}}
\oplus E_{j}^{(1)} \oplus \cdots \oplus E_{j}^{(n_j)})$,
where $\underline{\mathbb{C}}$ is the trivial line bundle over $B_{j-1}$,
$E_{j}^{(k)}$ is a holomorphic line bundle over $B_{j-1}$ for $k = 1, \dots, n_j$.
We call $B_j$ a
{\it $j$-stage generalized Bott manifold}. A generalized Bott tower has been
studied from various points of view (see, e.g., \cite{CMS,CMS2,HLS}).
A generalized Bott manifold is a certain class of toric manifold,
so it is interesting to investigate the properties of generalized Bott towers.

In~\cite{GK}, Grossberg and Karshon showed the multiplicity function of the
{\it equivariant index} (see $\S 2.4$) for a holomorphic line bundle over a Bott
manifold is given by the
density function of a {\it twisted cube}, which is determined by the structure
of the Bott manifold and the line bundle over it. From this, they derived a
Demazure-type character formula.

The purpose of this paper is to generalize the results in~\cite{GK} to
generalized Bott manifolds.
We generalize the twisted cube, and we call it the {\it generalized twisted cube}.
It is a special case of twisted polytope introduced by Karshon and Tolman \cite{KT}
for the presymplectic toric manifold, and it is a special case of
multi-polytope introduced by
Hattori and Masuda \cite{HM} for the torus manifold.
We show the multiplicity function of the equivariant index for the holomorphic line
bundle over the generalized Bott manifold is given by the density function of the
generalized twisted cube. From this, we derive a Demazure-type character formula.
In order to state the main results, we give some notation.
Let $\mathbf{L}$ be a holomorphic line bundle over a generalized Bott manifold
$B_m$, which is
constructed from integers $\{ \ell_i \}$ and $\{ c_{i,j}^{(k)} \}$ (see $\S 2.1$).
Let $N = \sum_{j=1}^m n_j$, and let $T^N = S^1 \times \cdots \times S^1$.
We consider the action of $T^N$ on $B_m$ as follows:
\begin{equation*}
(\mathbf{t}_1, \dots, \mathbf{t}_m) \cdot [\mathbf{z}_1, \dots, \mathbf{z}_m]
= [\mathbf{t}_{1}\mathbf{z}_1, \dots, \mathbf{t}_{m}\mathbf{z}_m],
\end{equation*}
where $\mathbf{t}_i = (t_{i,1}, \dots, t_{i,n_i}), \mathbf{z}_i = (z_{i,0}, \dots,
z_{i,n_i}), \mathbf{t}_{i}\mathbf{z}_i = (z_{i,0}, t_{i,1}z_{i,1}, \dots,
t_{i,n_i}z_{i,n_i})$ for $i = 1, \dots, m$.
Also we consider the action of $T = T^N \times S^1$ on $\mathbf{L}$ as follows:
\begin{equation}\label{(0.1)}
(\mathbf{t}_{1}, \dots, \mathbf{t}_{m}, t_{m+1}) \cdot
[\mathbf{z}_{1}, \dots, \mathbf{z}_{m}, v]
= [\mathbf{t}_{1}\mathbf{z}_{1}, \dots, \mathbf{t}_{m}\mathbf{z}_{m}, t_{m+1}v].
\end{equation}
We define the generalized twisted cube as follows. It is defined
to be the set of $x = (x_{1,1}, \dots, x_{m,n_m}) \in \mathbb{R}^N$ which satisfies
\begin{align*}
&A_i(x) \leq \sum_{k=1}^{n_i} x_{i,k} \leq 0, \,\, x_{i,k} \leq 0 \,\,\,\,
(1 \leq k \leq n_i) \notag \\
&{\rm or} \,\, 0 < \sum_{k=1}^{n_i} x_{i,k} < A_i(x), \,\, x_{i,k} > 0 \,\,\,\,
(1 \leq k \leq n_i),
\end{align*}
for $1 \leq i \leq m$, where
\begin{equation*}
A_i(x) =
\begin{cases}
-\ell_m & (i=m) \\
-(\ell_{i} + \sum_{j=i+1}^{m}
\sum_{k=1}^{n_j} c_{i,j}^{(k)}x_{j,k}) & (1 \leq i \leq m-1).
\end{cases}
\end{equation*}
We denote the generalized twisted cube by $C$.
We also define ${\rm sgn}(x_{i,k}) = 1$ for $x_{i,k} > 0$
and ${\rm sgn}(x_{i,k}) = -1$ for $x_{i,k} \leq 0$.
The {\it density function} of the generalized twisted cube is defined to be
$\rho(x) = (-1)^N \prod_{i,k} {\rm sgn}(x_{i,k})$ when
$x \in C$ and $0$ elsewhere.

Let $\mathfrak{t}$ be the Lie algebra of $T$ and let $\mathfrak{t}^{\ast}$ be its dual
space. Let $\ell^{\ast} \subset i\mathfrak{t}^{\ast}$ be the integral weight lattice
and let ${\rm mult} : \ell^{\ast} \to \mathbb{Z}$ be the multiplicity function of the
equivariant index.
The first main result of this paper is the following:
\begin{thm}\label{thm1.1}
Fix integers $\{ c_{i,j}^{(k)} \}$ and $\{ \ell_j \}$. Let $\mathbf{L} \to B_m$ be the
corresponding line bundle over a generalized Bott manifold.
Let $\rho \,:\, \mathbb{R}^N \to\{ -1, 0, 1 \}$
be the density function of the generalized twisted cube $C$ which is
determined by these integers.
Consider the torus action of
$T = T^N \times S^1$ as in (\ref{(0.1)}). Then the multiplicity function
for $\ell^{\ast} \cong \mathbb{Z}^N \times \mathbb{Z}$
is given by
\begin{equation*}
{\rm mult}(x,k)
= \begin{cases}
\rho(x) & (k = 1) \\
0 & (k \neq 1).
\end{cases}
\end{equation*}
\end{thm}
Karshon and Tolman found a toric manifold for which the multiplicities
of the equivariant index are $0,-1$, or $-2$ (\cite[Example 6.7]{KT}).
A generalized Bott manifold is different from this case by Theorem~\ref{thm1.1}.

Next, we give our character formula. Let
$\{ e_{1,1}, \dots, e_{m,n_m}, e_{m+1} \}$ be the standard basis in
$\mathbb{R}^{N+1}$, $x_i = (x_{i,1},
\dots, x_{i,n_i})$, and $e_i = (e_{i,1}, \dots, e_{i,n_i})$.
Let $\Delta_{n,r}^{-} = \left\{ z = (z_1, \dots, z_n) \in \mathbb{Z}_{\leq 0}^{n}
\,\middle|\, \right.$ \\
$\left. z_1 + \cdots + z_n = -r \right\}$,
and let
$\Delta_{n,r}^{+} = \left\{ z = (z_1, \dots, z_n) \in \mathbb{Z}_{> 0}^n
\,\middle|\, z_1 + \cdots + z_n = r-1 \right\}$.
For every integral weight
$\mu \in \ell^{\ast}$ we have a homomorphism $\lambda^{\mu} : T \to S^1$. We denote
the integral combinations of these $\lambda^{\mu}$'s by $\mathbb{Z}[T]$.
Then the operators $D_i : \mathbb{Z}[T] \to \mathbb{Z}[T]$ are defined using
$c_{i,j}^{(k)}$ and $\ell_j$
in the following way:
\begin{equation*}
D_i(\lambda^{\mu}) = \begin{cases}
\displaystyle \sum_{0 \leq r \leq k_i} \sum_{x_i \in \Delta_{n_i,r}^{-}}
\lambda^{\mu + \langle x_i, e_i \rangle} & \text{if} \,\, k_i \geq 0 \\
0 & \text{if} \,\, -n_i \leq k_i \leq -1 \\
\displaystyle \sum_{n_i + 1 \leq r \leq -k_i} \sum_{x_i \in \Delta_{n_i,r}^{+}}
(-1)^{n_i}\lambda^{\mu + \langle x_i, e_i \rangle} & \text{if} \,\, k_i \leq -n_i - 1,
\end{cases}
\end{equation*}
where the functions $k_i$ are defined as follows: if $\mu = e_{m+1} + \sum_{j=i+1}^m
\sum_{k=1}^{n_j} x_{j,k}e_{j,k}$,
then $k_i(\mu) = \ell_i + \sum_{j=i+1}^m \sum_{k=1}^{n_j} c_{i,j}^{(k)}x_{j,k}$.
From Theorem~\ref{thm1.1}, we obtain the following theorem:
\begin{thm}
Consider the action of the torus $T$ on $\mathbf{L} \to B_m$ as in (\ref{(0.1)}).
Denote the $(N+1)$-th component of the standard basis in $\mathbb{R}^{N+1}$ by $e_{m+1}$.
Then the character is given by the following element of $\mathbb{Z}[T]$:
\[ \chi = D_1 \cdots D_m(\lambda^{e_{m+1}}).
\]
\end{thm}
This is a Demazure-type character formula. On the other hand, the character is also
given by the localization formula with respect to the action of $T$
(\cite[Corollary 7.4]{HM}).
We compare our formula and the localization formula (see Remark~\ref{rmk3.7}).

This paper is organized as follows. In Section 2, we recall the equivariant index and
the generalized Bott towers, and we give the definition of
generalized twisted cubes. In Section 3, we prove the main theorems.


\section{\large Preliminaries}

In this section, we set up the tools to prove the main theorems.

\subsection{Generalized Bott manifolds}

\begin{defi}[\cite{CMS}]
\textup{A {\it generalized Bott tower} of height $m$ is a sequence:}
\[ B_m \stackrel{\pi_m}{\to} B_{m-1} \stackrel{\pi_{m-1}}{\to} \cdots \stackrel{\pi_2}{\to} B_1 \stackrel{\pi_1}{\to} B_0 = \{ \textup{a point} \},
\]
\textup{of manifolds $B_j = \mathbb{P}(\underline{\mathbb{C}} \oplus E_{j}^{(1)} \oplus \cdots \oplus E_{j}^{(n_j)})$,
where $\underline{\mathbb{C}}$ is the trivial line bundle over $B_{j-1}$, $E_{j}^{(k)}$ is a holomorphic line bundle over $B_{j-1}$ for $k = 1, \dots n_j$,
and $\mathbb{P}(\cdot)$ denotes the projectivization. We call $B_j$ a} $j$-stage generalized Bott manifold.
\end{defi}

The construction of the generalized Bott tower is as follows. A 1-step generalized
Bott tower can be written as
$B_1 = \mathbb{C}P^{n_1} = (\mathbb{C}^{n_1 + 1})^{\times} / \mathbb{C}^{\times}$,
where $\mathbb{C}^{\times}$ acts diagonally. We construct a line bundle over $B_1$ by
$E_{2}^{(k)} = (\mathbb{C}^{n_1 + 1})^{\times} \times_{\mathbb{C}^{\times}} \mathbb{C}$
for $k = 1, \dots, n_2$,
where $\mathbb{C}^{\times}$ acts on $\mathbb{C}$ by $a : v \mapsto a^{-c_k}v$
for some integer $c_k$.
In $E_{2}^{(k)}$ we have $[z_{1,0}, \dots, z_{1,n_1},v] = [z_{1,0}a, \dots,
z_{1,n_1}a, a^{c_k}v]$ for all $a \in \mathbb{C}^{\times}$.
A 2-step generalized Bott tower $B_2 = \mathbb{P}(\underline{\mathbb{C}} \oplus
E_{2}^{(1)} \oplus \cdots \oplus E_{2}^{(n_2)})$
can be written as $B_2 = ((\mathbb{C}^{n_1 + 1})^{\times} \times
(\mathbb{C}^{n_2 + 1})^{\times}) / G$,
where the right action of $G = (\mathbb{C}^{\times})^2$ is given by
\begin{equation*}
(\mathbf{z}_1,\mathbf{z}_2) \cdot (a_1,a_2) =
(z_{1,0}a_1, z_{1,1}a_1, \dots, z_{1,n_1}a_1,z_{2,0}a_2, a_{1}^{c_1}z_{2,1}a_2,
\dots, a_{1}^{c_{n_2}}z_{2,n_2}a_2),
\end{equation*}
where $\mathbf{z}_j = (z_{j,0}, z_{j,1}, \dots, z_{j,n_j})$ for $j = 1,2$.

We can construct higher generalized Bott tower in a similar way.
In this way we get an $m$-step generalized Bott manifold
$B_m = \mathbb{P}(\underline{\mathbb{C}} \oplus E_{m}^{(1)} \oplus \cdots \oplus
E_{m}^{(n_m)})$
from any collection of integers $\{ c_{i,j}^{(k)} \}$:
\[ B_m = ((\mathbb{C}^{n_1 + 1})^{\times} \times \cdots \times (\mathbb{C}^{n_m + 1})^{\times}) / G,
\]
where the right action of $G = (\mathbb{C}^{\times})^m$ is given by
\begin{equation*}
(\mathbf{z}_1, \dots, \mathbf{z}_m) \cdot \mathbf{a}
= (\mathbf{z}_{1}', \mathbf{z}_{2}', \dots, \mathbf{z}_{m}'),
\end{equation*}
where $\mathbf{z}_i = (z_{i,0}, \dots, z_{i,n_i})$ for $i = 1, \dots, m$,
$\mathbf{a} = (a_1, \dots, a_m) \in (\mathbb{C}^{\times})^m$, \\
$\mathbf{z}_1' = (z_{1,0}a_1, z_{1,1}a_1, \dots, z_{1,n_1}a_1)$
and $\mathbf{z}_j' = (z_{j,0}a_j,
a_{1}^{c_{1,j}^{(1)}} \cdots a_{j-1}^{c_{j-1,j}^{(1)}}z_{j,1}a_j, \dots,
a_{1}^{c_{1,j}^{(n_j)}} \cdots a_{j-1}^{c_{j-1,j}^{(n_j)}}z_{j,n_j}a_j)$ for
$j = 2, \dots, m$.
We can construct a line bundle over $B_m$ from the integers $(\ell_1, \dots, \ell_m)$ by
\[ \mathbf{L} = ((\mathbb{C}^{n_1 + 1})^{\times} \times \cdots \times (\mathbb{C}^{n_m + 1})^{\times}) \times_{G} \mathbb{C},
\]
where $G = (\mathbb{C}^{\times})^m$ acts by
\begin{equation}\label{(1)}
((\mathbf{z}_1, \dots, \mathbf{z}_m), v) \cdot \mathbf{a}
= (\mathbf{z}_{1}', \mathbf{z}_{2}', \dots, \mathbf{z}_{m}', a_{1}^{\ell_1} \cdots
a_{m}^{\ell_m}v).
\end{equation}


\subsection{Torus action on generalized Bott towers}

Let $N = \sum_{j=1}^{m} n_j$ and let $T^N = S^1 \times \cdots \times S^1$.
We consider the action of $T^N$ on $B_m$ as follows:
\begin{equation*}
(\mathbf{t}_1, \dots, \mathbf{t}_m) \cdot [\mathbf{z}_1, \dots, \mathbf{z}_m]
= [\mathbf{t}_{1} \cdot \mathbf{z}_1, \dots, \mathbf{t}_m \cdot \mathbf{z}_m],
\end{equation*}
where $\mathbf{t}_i = (t_{i,1}, \dots, t_{i,n_i})$ and
$\mathbf{t}_i \cdot \mathbf{z}_i = (z_{i,0}, t_{i,1}z_{i,1}, \dots,
t_{i,n_i}z_{i,n_i})$ for $i = 1, \dots, m$.
Also we consider the action of $T = T^N \times S^1$ on $\mathbf{L}$ as follows:
\begin{equation}\label{(2.1)}
(\mathbf{t}_{1}, \dots, \mathbf{t}_{m}, t_{m+1}) \cdot [\mathbf{z}_{1}, \dots,
\mathbf{z}_{m}, v]
= [\mathbf{t}_{1} \cdot \mathbf{z}_1, \dots, \mathbf{t}_m \cdot \mathbf{z}_m,
t_{m+1}v].
\end{equation}


\subsection{Generalized twisted cubes}

\begin{defi}
\textup{A {\it generalized twisted cube} $C$ is defined to be the
set of $x = (x_{1,1}, \dots, x_{m,n_m})$ \\
$\in \mathbb{R}^N$ which satisfies}
\begin{align}\label{(3.1)}
&A_i(x) \leq \sum_{k=1}^{n_i} x_{i,k} \leq 0, \,\,
x_{i,k} \leq 0 \,\,\,\, (1 \leq k \leq n_i) \notag \\
&{\rm or} \,\, 0 < \sum_{k=1}^{n_i} x_{i,k} < A_i(x),
\,\, x_{i,k} > 0 \,\,\,\, (1 \leq k \leq n_i),
\end{align}
\textup{for all $1 \leq i \leq m$, where}
\begin{equation*}
A_i(x) =
\begin{cases}
-\ell_m & (i=m) \\
-(\ell_{i} + \sum_{j=i+1}^{m} \sum_{k=1}^{n_j} c_{i,j}^{(k)}x_{j,k}) &
(1 \leq i \leq m-1).
\end{cases}
\end{equation*}
\end{defi}

\begin{rema}
\textup{(i) The generalized twisted cube is a special case of multi-polytope defined in \cite{HM}.
In particular, it is a special case of twisted polytope defined in \cite{KT}.}

\textup{(ii) When $n_i = 1$ for all $1 \leq i \leq m$, the generalized twisted cube
is the twisted cube given in~\cite[(2.21)]{GK}.}
\end{rema}

\begin{defi}
\textup{We define ${\rm sgn}(x_{i,k}) = 1$ for $x_{i,k} > 0$ and ${\rm sgn}(x_{i,k}) = -1$ for $x_{i,k} \leq 0$.
The {\it density function} of the generalized twisted cube
is then defined to be $\rho(x) = (-1)^N \prod_{i,k} {\rm sgn}(x_{i,k})$
when $x \in C$ and $0$ elsewhere.}
\end{defi}

\begin{exa}\label{ex2}
{\rm Suppose that $m=2, n_1 = 1, n_2 = 2, \ell_1 = 1$, and $\ell_2 = 2$. We set $c_{1,2}^{(1)} = 2$ and $c_{1,2}^{(2)} = -1$.
Then the generalized twisted cube is the set of
 $x = (x_{1,1}, x_{2,1}, x_{2,2})$ which satisfies}
\begin{itemize}
\item $-2 \leq x_{2,1} + x_{2,2} \leq 0, \,\, x_{2,1}, x_{2,2} \leq 0,$
\item $-1 - 2x_{2,1} + x_{2,2} \leq x_{1,1} \leq 0 \,\,
\textup{or} \,\, 0 < x_{1,1} < -1 - 2x_{2,1} + x_{2,2}.$
\end{itemize}
{\rm In Figure 1, the black dots represent the lattice points of the sign $+1$ and the
white dots represent the sign $-1$.}

\begin{center}
\begin{tikzpicture}
\draw[line width=1.5pt] (0,-2) -- (0,0) -- (-2,0) -- (-1,-1) ;
\draw[line width=1.5pt] (0,0) -- (-0.5,-0.8) -- (-1.5,-4.4) -- (0,-2) -- cycle ;
\draw[line width=1.5pt] (-0.5,-0.8) -- (-0.6,1.6) -- (-2,0) ;
\draw[line width=1.5pt] (-0.5,0) -- (-1,-1) ;
\draw[dashed] (-0.6,1.6) -- (-1,-1) ;
\draw[line width=1.5pt] (-1,-1) -- (-1.5,-4.4) ;
\draw[dashed] (-1,-1) -- (0,-2) ;
\fill (0,-1) circle (3pt) ;
\fill (0,0) circle (3pt) ;
\fill (0,-2) circle (3pt) ;
\fill (-1.5,-4.4) circle (3pt) ;
\fill (-0.5,-0.8) circle (3pt) ;
\fill (-1,-2.6) circle (3pt) ;
\fill (-0.5,-1.8) circle (3pt) ;
\fill (-0.5,-2.8) circle (3pt) ;
\fill (-1,-3.6) circle (3pt) ;
\fill (-1,-1) circle (3pt) ;
\draw (-1.5,0.6) circle (3pt) ;
\draw (-1,1.1) circle (3pt) ;
\draw (0.4,0) node {$\textup{O}$} ;
\draw (1,-2) node {$(0,0,-2)$} ;
\draw (-0.3,-4.4) node {$(-3,0,-2)$} ;
\draw (-4.5,0) node {$(x_{1,1},x_{2,1},x_{2,2}) = (0,-2,0)$} ;
\draw (0.3,1.6) node {$(3,-2,0)$} ;
\end{tikzpicture}
\end{center}
\begin{center}
{\rm Figure 1}
\end{center}
\end{exa}

\begin{exa}\label{ex3}
{\rm Suppose that $m=2, n_1 = 2, n_2 = 1, \ell_1 = 2$, and $\ell_2 = -6$.
We set $c_{1,2}^{(1)} = -1$. Then the generalized twisted cube is the set of
$x = (x_{1,1}, x_{1,2}, x_{2,1})$ which satisfies}
\begin{itemize}
\item $0 < x_{2,1} < 6,$
\item $-2 + x_{2,1} \leq x_{1,1} + x_{1,2} \leq 0, \,\, x_{1,1}, x_{1,2} \leq 0 \,\,
\textup{or} \,\, 0 < x_{1,1} + x_{1,2} < -2 + x_{2,1}, \,\, x_{1,1}, x_{1,2} > 0.$
\end{itemize}
{\rm In Figure 2, the white dots represent the sign $-1$.}

\begin{center}
\begin{tikzpicture}
\draw[line width=1.5pt] (0,0) -- (-2,0) -- (0,2) -- (0,0) ;
\draw[line width=1.5pt] (0,6) -- (4,6) -- (0,2) ;
\draw[dashed] (0,2) -- (0,6) ;
\draw[line width=1.5pt] (0,0) -- (0.8,1.2) -- (0,2) ;
\draw[dashed] (0.8,1.2) -- (-2,0) ;
\draw[line width=1.5pt] (0,6) -- (-1.6,3.6) -- (0,2) ;
\draw[line width=1.5pt] (-1.6,3.6) -- (4,6) ;
\draw (-1,1) circle (3pt) ;
\draw (0,1) circle (3pt) ;
\draw (0.4,1.6) circle (3pt) ;
\draw (0,2) circle (3pt) ;
\draw (0.6,4.4) circle (3pt) ;
\draw (0.3,0) node {$\textup{O}$} ;
\draw (-3,0) node {$(0,-2,0)$} ;
\draw (-1,2) node {$(0,0,2)$} ;
\draw (-2.5,6) node {$(x_{1,1},x_{1,2},x_{2,1}) = (0,0,6)$} ;
\draw (5,6) node {$(0,4,6)$} ;
\draw (-2.6,3.6) node {$(4,0,6)$} ;
\draw (1.5,4.4) node {$(1,1,5)$} ;
\end{tikzpicture}
\end{center}
\begin{center}
{\rm Figure 2}
\end{center}
\end{exa}

\begin{exa}\label{ex4}
{\rm Suppose that $m=2,n_1=n_2=2,\ell_1=1$, and $\ell_2=2$. We set
$c_{1,2}^{(1)}=2$ and $c_{1,2}^{(2)}=-1$. Then the generalized twisted cube is
the set of $x = (x_{1,1}, x_{1,2}, x_{2,1}, x_{2,2})$ which satisfies}
\begin{itemize}
\item $-2 \leq x_{2,1} + x_{2,2} \leq 0, \,\, x_{2,1}, x_{2,2} \leq 0,$
\item $-1-2x_{2,1}+x_{2,2} \leq x_{1,1} + x_{1,2} \leq 0, \,\, x_{1,1}, x_{1,2} \leq 0 \\
\textup{or} \,\, 0 < x_{1,1} + x_{1,2} < -1-2x_{2,1}+x_{2,2}, \,\,
x_{1,1}, x_{1,2} > 0.$
\end{itemize}
{\rm The lattice points in the generalized twisted cube represent the sign $-1$.}
\end{exa}


\subsection{Equivariant index}

Let $\mathbf{L}$ be a holomorphic line bundle over a generalized Bott manifold $B_m$
with the action of the torus $T$ as in (\ref{(2.1)}).
Let $\mathcal{O}_{\mathbf{L}}$ be the sheaf of holomorphic sections.
The {\it equivariant index} of a generalized Bott manifold
is the formal sum of representation of $T$:
\[
{\rm index}(B_m, \mathcal{O}_{\mathbf{L}}) = \sum (-1)^{i}H^i(B_m, \mathcal{O}_{\mathbf{L}})
\]
The {\it character} of the equivariant index is the function $\chi : T \to \mathbb{C}$
which is given by $\chi = \sum (-1)^i \chi^i$ where
$\chi^i(a) = \text{trace}\{ a : H^i(B_m, \mathcal{O}_{\mathbf{L}}) \to
H^i(B_m, \mathcal{O}_{\mathbf{L}}) \}$ for $a \in T$.
Let $\mathfrak{t}$ be the Lie algebra of $T$ and let $\mathfrak{t}^{\ast}$ be its
dual space.
Every $\mu$ in the integral
weight lattice $\ell^{\ast} \subset i\mathfrak{t}^{\ast}$ defines a homomorphism
$\lambda^{\mu} : T \to S^1$. We can write
$\chi = \sum_{\mu \in \ell^{\ast}} m_{\mu}\lambda^{\mu}$. The coefficients
are given by a function $\text{mult} \,: \ell^{\ast} \to \mathbb{Z}$,
sending $\mu \mapsto m_{\mu}$, called the {\it multiplicity function}
for the equivariant index.


\section{\large Main theorems}

\subsection{Multiplicity function of the equivariant index}

We will show that the multiplicity function of the
equivariant index of a generalized Bott manifold is given by the density function of a
generalized twisted cube $C$.
In particular, all the weights occur with a multiplicity $-1, 0$, or $1$.

\begin{thm}\label{thm1}
Fix integers $\{ c_{i,j}^{(k)} \}$ and $\{ \ell_j \}$. Let $\mathbf{L} \to B_m$ be
the corresponding line bundle over a generalized Bott manifold.
Let $\rho \,:\, \mathbb{R}^N \to \{ -1, 0, 1 \}$ be the density function of
the generalized twisted cube $C$ which is determined by these integers as in (\ref{(3.1)}).
Consider the torus action of
$T = T^N \times S^1$ as in (\ref{(2.1)}).
Then the multiplicity function for
$\ell^{\ast} \cong \mathbb{Z}^N \times \mathbb{Z}$ is given by
\begin{equation*}
{\rm mult}(x,k)
= \begin{cases}
\rho(x) & (k = 1) \\
0 & (k \neq 1).
\end{cases}
\end{equation*}
\end{thm}
\noindent
\textit{Proof} ;
We compute $H^{\ast}(B_m, \mathcal{O}_{\mathbf{L}})$.
Take the covering $\tilde{\mathcal{U}} = \{ U_{r_1} \times \cdots \times U_{r_m} \}$ of
$(\mathbb{C}^{n_1 + 1})^{\times} \times \cdots \times (\mathbb{C}^{n_m + 1})^{\times}$ for
$r_1, \dots, r_m \in \{ 0, 1, \dots, n_{\ell} \} \,\, (\ell = 1, \dots, m)$, where
$U_{r_j} = \underbrace{\mathbb{C} \times \cdots \times \mathbb{C}}_{r_j} \times
\mathbb{C}^{\times} \times \underbrace{\mathbb{C} \times \cdots \times \mathbb{C}}_{n_{\ell}-r_j}$.
This descends to the covering $\mathcal{U}$ of $B_m$; every intersection of sets in $\mathcal{U}$ is isomorphic to a
product of $\mathbb{C}$'s and $\mathbb{C}^{\times}$'s.
The coverings $\tilde{\mathcal{U}}$ and $\mathcal{U}$ are the Leray coverings
(\cite{GH}).

Let $\mathcal{O}$ be the sheaf of holomorphic functions, and let
$G = (\mathbb{C}^{\times})^m$. Since holomorphic sections of
$\mathcal{O}_{\mathbf{L}}$ are given by holomorphic sections of
$\mathcal{O}$ which are $G$-invariant with respect to the action (\ref{(1)})
(\cite{KT}), $H^{\ast}(\mathcal{U},\mathcal{O}_{\mathbf{L}})$
is isomorphic to the $G$-invariant part of $H^{\ast}(\tilde{\mathcal{U}},\mathcal{O})$.
By the Leray theorem, $H^{\ast}(B_m, \mathcal{O}_{\mathbf{L}})$
is isomorphic to the $G$-invariant part of
$H^{\ast}((\mathbb{C}^{n_1 + 1})^{\times} \times \cdots \times
(\mathbb{C}^{n_m + 1})^{\times}, \mathcal{O})$.

In order to compute $H^{\ast}((\mathbb{C}^{n_1 + 1})^{\times} \times \cdots \times
(\mathbb{C}^{n_m + 1})^{\times}, \mathcal{O})$,
we compute $H^{\ast}((\mathbb{C}^{n+1})^{\times}, \mathcal{O})$.
Let $\mathcal{U}' = \{ U_0,U_1, \dots, U_n \}$ be the covering of
$(\mathbb{C}^{n+1})^{\times}$, let $j_0, j_1, \dots, j_k \in \{ 0, 1, \dots, n \}$ for
$k = 0, 1, \dots, n$ and let $U_{j_{0}j_{1} \cdots j_{k}}
= U_{j_0} \cap U_{j_1} \cap \cdots \cap U_{j_k}$.
Let $I = (i_0, i_1, \dots, i_n) \in \mathbb{Z}^{n+1}$.
The holomorphic functions on $U_{j_{0}j_{1} \cdots j_{k}}$
are given by
\[
\Gamma_{\rm hol}(U_{j_{0}j_{1} \cdots j_{k}}) = \left\{ \sum_{I \in \mathbb{Z}^{n+1}, i_{\ell} \geq 0 (\ell \neq j_0, j_1, \dots, j_k)}
a_{I}z_{0}^{i_0}z_{1}^{i_1} \cdots z_{n}^{i_n} \right\}.
\]
Consider the \v{C}ech cochain complex
\[ 0 \to \check{C}^0(\mathcal{U}', \mathcal{O}) \stackrel{\delta^0}{\to} \check{C}^1(\mathcal{U}', \mathcal{O}) \stackrel{\delta^1}{\to} \cdots \stackrel{\delta^{n-1}}{\to}
\check{C}^n(\mathcal{U}', \mathcal{O}) \stackrel{\delta^n}{\to} 0,
\]
where $\check{C}^i(\mathcal{U}', \mathcal{O})
= \oplus \Gamma_{\rm hol}(U_{j_{0}j_{1} \cdots j_i}) \,\, (i = 0, \dots, n)$.
The map $\delta^p : \check{C}^p(\mathcal{U}', \mathcal{O}) \to
\check{C}^{p+1}(\mathcal{U}', \mathcal{O})$
is given by
$\{ f_{j_{0}j_{1} \cdots j_{p}}
\} \mapsto \{ g_{j_{0}j_{1} \cdots j_{p+1}} \}, g_{j_{0}j_{1} \cdots j_{p+1}}
= \sum (-1)^k f_{j_{0}j_{1} \cdots \hat{j_{k}} \cdots j_{p+1}}$.
Recall that \\
$H^0((\mathbb{C}^{n+1})^{\times}, \mathcal{O}) = {\rm Ker} \, \delta^0$,
 and $H^n((\mathbb{C}^{n+1})^{\times}, \mathcal{O})
= {\rm Coker} \, \delta^{n-1}$.
The torus $T^{n+1} = (S^1)^{n+1}$ acts on the holomorphic functions by
$((t_0, \dots, t_n) \cdot f)(z_0, \dots, z_n) = f(t_{0}^{-1}z_0, \dots, t_{n}^{-1}z_n)$.
This action descends to the cohomology. The corresponding weight spaces for the weight
$I \in \mathbb{Z}^{n+1}$ are
\begin{align*}
H^0((\mathbb{C}^{n+1})^{\times}, \mathcal{O})_{I} &=
\begin{cases}
{\rm span}(z_{0}^{-i_0} \cdots z_{n}^{-i_n}) &(I \in \mathbb{Z}_{\leq 0}^{n+1}) \\
0 &\text{otherwise}
\end{cases} \\
H^n((\mathbb{C}^{n+1})^{\times}, \mathcal{O})_I &=
\begin{cases}
{\rm span}(z_{0}^{-i_0} \cdots z_{n}^{-i_n}) &(I \in \mathbb{Z}_{>0}^{n+1}) \\
0 &\text{otherwise}.
\end{cases}
\end{align*}
We now prove $H^q((\mathbb{C}^{n+1})^{\times},\mathcal{O}) = 0$ for $1 \leq q \leq n-1$.
Let $\Delta$ be the fan of $(\mathbb{C}^{n+1})^{\times}$, and let $|\Delta| =
\cup_{\sigma \in \Delta} \sigma$.
Let
\begin{equation*}
Z(I) := \{ v \in |\Delta| \,;\, \langle I,v \rangle \leq \varphi(v) \},
\end{equation*}
where $\varphi$ is the support function. From~\cite{F},
\begin{equation*}
H^q((\mathbb{C}^{n+1})^{\times},\mathcal{O})_I
= H^q(|\Delta|,|\Delta| \setminus Z(I) \,;\, \mathbb{C}).
\end{equation*}
Since $\mathcal{O}$ is the sheaf of holomorphic function, $\varphi(v) = 0$ for all
$v \in |\Delta|$.
In the case that $i_j \leq 0$ for all $j$, since $|\Delta|$ is contractible,
\begin{equation*}
H^q((\mathbb{C}^{n+1})^{\times},\mathcal{O})_I = 0 \,\, (q \geq 1).
\end{equation*}
In the case that $i_j > 0$ for all $j$, $Z(I) = \{ 0 \}$. Since $|\Delta| \setminus
\{ 0 \}$ is homotopic to $S^{n-1}$,
\begin{equation*}
H^q((\mathbb{C}^{n+1})^{\times},\mathcal{O})_I = 0 \,\, (q \neq n).
\end{equation*}
In other case, since $|\Delta| \setminus Z(I)$ is path-connected and contractible,
\begin{equation*}
H^q((\mathbb{C}^{n+1})^{\times},\mathcal{O})_I = 0
\end{equation*}
for all $q$.

We now compute $H^{\ast}((\mathbb{C}^{n_1 + 1})^{\times} \times \cdots \times
(\mathbb{C}^{n_m + 1})^{\times}, \mathcal{O})$.
Consider the natural action of $T^{N+m} = (S^1)^{N+m}$ on the holomorphic function.
The weights are multi-indices $I' \in \mathbb{Z}^{N+m}$; we write
$I' = (\mathbf{i}'_{1}, \dots, \mathbf{i}'_{m})$,
where $\mathbf{i}'_{j} = (i_{j,0}, i_{j,1}, \dots, i_{j,n_j})$ for $j = 1, \dots, m$.
From the cohomology of $(\mathbb{C}^{n+1})^{\times}$ that we have computed and from the
K\"{u}nneth formula (\cite{AN}), it follows that
\begin{equation*}
H^{q}((\mathbb{C}^{n_1 + 1})^{\times} \times \cdots \times (\mathbb{C}^{n_m +
1})^{\times}, \mathcal{O})_{I'} =
\begin{cases}
\text{span}(z_{1,0}^{-i_{1,0}}z_{1,1}^{-i_{1,1}} \cdots z_{m,n_m}^{-i_{m,n_m}}) \\
0.
\end{cases}
\end{equation*}
The former occurs if for all $\ell$ we have
${\rm sgn}(i_{\ell,0}) = {\rm sgn}(i_{\ell,1}) = \cdots = {\rm sgn}(i_{\ell,n_{\ell}})
=: \varepsilon_{\ell}$, here $q = \sum_{\{ \ell \,|\, \varepsilon_{\ell} = 1, 1 \leq
\ell \leq m \}} n_{\ell}$,
and $q = 0$ when $\varepsilon_{\ell} = -1$ for all $\ell$. In particular, $(-1)^q = (-
1)^N \prod_{1 \leq \ell \leq m, 1 \leq p \leq n_{\ell}} {\rm sgn}(i_{\ell,p})$.

The action (\ref{(1)}) induces an action on functions given by
\begin{align*}
(a_{k}f)(z_{1,0}, \dots, z_{m,n_m}) = a_{k}^{\ell_k}f(z_{1,0}, \dots, z_{k-1,n_{k-1}},
&z_{k,0}a_{k}^{-1}, z_{k,1}a_{k}^{-1}, \dots, z_{k,n_k}a_{k}^{-1}, \\
&\dots, z_{\ell,0}, a_{k}^{-c_{k,\ell}^{(1)}}z_{\ell,1}, \dots, a_{k}^{-
c_{k,\ell}^{(n_{\ell})}}z_{\ell,n_{\ell}}, \dots).
\end{align*}
The monomial $z_{1,0}^{-i_{1,0}}z_{1,1}^{-i_{1,1}} \cdots z_{m,n_m}^{-i_{m,n_m}}$ is
then a weight vector with a weight whose
$k$-th coordinate is $\ell_k + i_{k,0} + \cdots + i_{k,n_k} + \sum_{\ell=k+1}^{m}
\sum_{p=1}^{n_{\ell}} c_{k,\ell}^{(p)}i_{\ell,p}$. Thus the $G$-invariant part of
$H^{\ast}((\mathbb{C}^{n_1 + 1})^{\times} \times \cdots \times (\mathbb{C}^{n_m +
1})^{\times}, \mathcal{O})$ consists of those monomials $z_{1,0}^{-i_{1,0}}z_{1,1}^{-
i_{1,1}} \cdots z_{m,n_m}^{-i_{m,n_m}}$ for which
\begin{align}\label{(2)}
\ell_1 + i_{1,0} + \cdots + i_{1,n_1} + \sum_{\ell=2}^{m} \sum_{p=1}^{n_{\ell}}
c_{1,\ell}^{(p)}i_{\ell,p} &= 0 \notag \\
\ell_2 + i_{2,0} + \cdots + i_{2,n_2} + \sum_{\ell=3}^{m} \sum_{p=1}^{n_{\ell}}
c_{2,\ell}^{(p)}i_{\ell,p} &= 0 \\
&\vdots \notag \\
\ell_m + i_{m,0} + \cdots + i_{m,n_m} &= 0. \notag
\end{align}

The action (\ref{(2.1)}) induces a $T$ action on the functions given by
\begin{align*}
&((t_{1,1}, \dots, t_{m,n_m}, t_{m+1}) \cdot f)(z_{1,0}, \dots, z_{m,n_m}) \\
&= t_{m+1}f(z_{1,0}, t_{1,1}^{-1}z_{1,1}, \dots, z_{m,0}, t_{m,1}^{-1}z_{m,1}, \dots,
t_{m,n_m}^{-1}z_{m,n_m}).
\end{align*}
The weight of the monomial $z_{1,0}^{-i_{1,0}}z_{1,1}^{-i_{1,1}} \cdots z_{m,n_m}^{-
i_{m,n_m}}$ with respect to this $T$ action is \\
$(\mathbf{i}_{1}, \mathbf{i}_{2}, \dots, \mathbf{i}_{m}, 1)$,
where $\mathbf{i}_j = (i_{j,1}, \dots, i_{j,n_j})$ for $j = 1, \dots, m$.
Thus the index of $(B_m, \mathcal{O}_{\mathbf{L}})$ is given by the set of
$x = (x_{1,1}, \dots, x_{m,n_m}, 1) = (i_{1,1}, \dots, i_{m,n_m}, 1)$ for which
there exist $(i_{1,0}, \dots, i_{m,0})$ such that (\ref{(2)}) is satisfied and such that
${\rm sgn}(i_{\ell,0}) = {\rm sgn}(i_{\ell,1}) = \cdots =
{\rm sgn}(i_{\ell,n_{\ell}})$ for all $\ell$. This is exactly the set (\ref{(3.1)}).
Therefore the multiplicity of the equivariant index is $(-1)^N \prod_{1 \leq \ell \leq
m, 1 \leq p \leq n_{\ell}} {\rm sgn}(i_{\ell,p})
= (-1)^N \prod_{1 \leq \ell \leq m, 1 \leq p \leq n_{\ell}} {\rm sgn}(x_{\ell,p}) =
\rho(x)$. \,\, $\square$

\subsection{Character formula for the equivariant index}

In the following the theorem we give a formula for the character
$\chi \,:\, T \to \mathbb{C}$ of the equivariant index of a
generalized Bott manifold.
For every integral weight $\mu \in \ell^{\ast}$ we have a homomorphism
$\lambda^{\mu} \,:\, T \to S^1$. We denote
the integral combinations of these $\lambda^{\mu}$'s by $\mathbb{Z}[T]$.
Then $\chi \in \mathbb{Z}[T]$ is given by $\chi = \sum_{\mu \in
\ell^{\ast}} m_{\mu}\lambda^{\mu}$
where $m_{\mu} = {\rm mult}(\mu)$.

\begin{defi}
{\rm Let $\{ e_{1,1}, \dots, e_{m,n_m}, e_{m+1} \}$ be the standard basis in
$\mathbb{R}^{N+1}$, $x_i = \left(x_{i,1}, \dots, \right.$ \\
$\left. x_{i,n_i} \right)$ and $e_i = (e_{i,1}, \dots, e_{i,n_i})$.
Let $\Delta_{n,r}^{-} = \left\{ z = (z_1, \dots, z_n) \in
\mathbb{Z}_{\leq 0}^{n} \,\middle|\, z_1
+ \cdots + z_n = -r \right\}$, and let
$\Delta_{n,r}^{+} = \left\{ z = (z_1, \dots, z_n) \in \mathbb{Z}_{> 0}^n \,\middle|\,
z_1 + \cdots +
z_n = r-1 \right\}$. Then the operators $D_i : \mathbb{Z}[T] \to \mathbb{Z}[T]$ are
defined using $c_{i,j}^{(k)}$ and $\ell_j$ in the following way:}
\begin{equation*}
D_i(\lambda^{\mu}) = \begin{cases}
\displaystyle \sum_{0 \leq r \leq k_i} \sum_{x_i \in \Delta_{n_i,r}^{-}}
\lambda^{\mu + \langle x_i, e_i \rangle} & \text{if} \,\, k_i \geq 0 \\
0 & \text{if} \,\, -n_i \leq k_i \leq -1 \\
\displaystyle \sum_{n_i + 1 \leq r \leq -k_i} \sum_{x_i \in \Delta_{n_i,r}^{+}}
(-1)^{n_i}\lambda^{\mu + \langle x_i, e_i \rangle} & \text{if} \,\, k_i \leq -n_i - 1,
\end{cases}
\end{equation*}
{\rm where the functions $k_i$ are defined as follows: if $\mu = e_{m+1} + \sum_{j=i+1}^m
\sum_{k=1}^{n_j} x_{j,k}e_{j,k}$, then $k_i(\mu) = \ell_i + \sum_{j=i+1}^m
\sum_{k=1}^{n_j} c_{i,j}^{(k)}x_{j,k}$.}
\end{defi}

From Theorem~\ref{thm1}, we immediately obtain the following theorem.

\begin{thm}\label{prop2}
Consider the action of the torus $T$ on $\mathbf{L} \to B_m$ as in (\ref{(2.1)}).
Denote the $(N+1)$-th component of the standard basis in $\mathbb{R}^{N+1}$ by
$e_{m+1}$.
Then the character is given by the following element of $\mathbb{Z}[T]$:
\[ \chi = D_1 \cdots D_m(\lambda^{e_{m+1}}).
\]
\end{thm}

\begin{rema}
{\rm When $n_i = 1$ for all $i$, the operator $D_i$ is given by}
\begin{equation*}
D_i(\lambda^{\mu}) = \begin{cases}
\lambda^{\mu} + \lambda^{\mu - e_{i,1}} + \cdots + \lambda^{\mu - k_{i}e_{i,1}} &
{\it if} \,\, k_i \geq 0 \\
0 & {\it if} \,\, k_i = -1 \\
-\lambda^{\mu + e_{i,1}} - \lambda^{\mu + 2e_{i,1}} - \cdots - \lambda^{\mu -
(k_i + 1)e_{i,1}} & {\it if} \,\, k_i \leq -2.
\end{cases}
\end{equation*}
{\rm We can check that this operator agrees with the one
in~\cite[Proposition 2.32]{GK}.}
\end{rema}

\begin{exa}\label{ex3.4}
{\rm Suppose that $m=2,n_1=1$, and $n_2=2$. We set $\ell_1=1,\ell_2=2,c_{1,2}^{(1)}=2$,
and $c_{1,2}^{(2)}=-1$ as in Example~\ref{ex2}. Then the corresponding
character $\chi$ is given by}
\begin{align*}
\chi &= D_1D_2(\lambda^{e_3}) \\
&= D_1(\lambda^{e_3} + \lambda^{e_3 - e_{2,1}} + \lambda^{e_3 - e_{2,2}} +
\lambda^{e_3 - 2e_{2,1}} + \lambda^{e_3 - e_{2,1} - e_{2,2}} +
\lambda^{e_3 - 2e_{2,2}}) \\
&= \lambda^{e_3} + \lambda^{e_3-e_{1,1}} + \lambda^{e_3-e_{2,2}}
+ \lambda^{e_3-e_{2,2}-e_{1,1}} + \lambda^{e_3-e_{2,2}-2e_{1,1}}
- \lambda^{e_3-2e_{2,1}+e_{1,1}} - \lambda^{e_3-2e_{2,1}+2e_{1,1}} \\
&\,\,\,\, + \lambda^{e_3-e_{2,1}-e_{2,2}} + \lambda^{e_3-2e_{2,2}}
+ \lambda^{e_3-2e_{2,2}-e_{1,1}} + \lambda^{e_3-2e_{2,2}-2e_{1,1}}
+ \lambda^{e_3-2e_{2,2}-3e_{1,1}}.
\end{align*}
\end{exa}

\begin{exa}
{\rm Suppose that $m=2,n_1=2$, and $n_2=1$. We set $\ell_1=2,\ell_2=-6$, and
$c_{1,2}^{(1)}=-1$ as in Example~\ref{ex3}. Then the corresponding
character $\chi$ is given by}
\begin{align*}
\chi &= D_1D_2(\lambda^{e_3}) \\
&= D_1(-\lambda^{e_3 + e_{2,1}} - \lambda^{e_3 + 2e_{2,1}} - \lambda^{e_3 + 3e_{2,1}}
- \lambda^{e_3 + 4e_{2,1}} - \lambda^{e_3 + 5e_{2,1}}) \\
&= -\lambda^{e_3+e_{2,1}} - \lambda^{e_3+e_{2,1}-e_{1,1}}
- \lambda^{e_3+e_{2,1}-e_{1,2}} - \lambda^{e_3+2e_{2,1}}
- \lambda^{e_3+5e_{2,1}+e_{1,1}+e_{1,2}}.
\end{align*}
\end{exa}

\begin{exa}
{\rm Suppose that $m=2,n_1=2$, and $n_2=2$. We set $\ell_1=1,\ell_2=2,
c_{1,2}^{(1)}=2$, and $c_{1,2}^{(2)}=-1$ as in Example~\ref{ex4}.
Then the corresponding character $\chi$ is given by}
\begin{align*}
\chi &= D_1D_2(\lambda^{e_3}) \\
&= D_1(\lambda^{e_3} + \lambda^{e_3-e_{2,1}} + \lambda^{e_3-e_{2,2}}
+ \lambda^{e_3-2e_{2,1}} + \lambda^{e_3-e_{2,1}-e_{2,2}} + \lambda^{e_3-2e_{2,2}}) \\
&= \lambda^{e_3} + \lambda^{e_3-e_{1,1}} + \lambda^{e_3-e_{1,2}}
+ \lambda^{e_3-e_{2,2}} + \lambda^{e_3-e_{2,2}-e_{1,1}} + \lambda^{e_3-e_{2,2}-e_{1,2}}
+ \lambda^{e_3-e_{2,2}-2e_{1,1}} \\
&\,\,\,\, + \lambda^{e_3-e_{2,2}-e_{1,1}-e_{1,2}}
+ \lambda^{e_3-e_{2,2}-2e_{1,2}} + \lambda^{e_3-2e_{2,1}+e_{1,1}+e_{1,2}}
+ \lambda^{e_3-e_{2,1}-e_{2,2}} + \lambda^{e_3-2e_{2,2}} \\
&\,\,\,\, + \lambda^{e_3-2e_{2,2}-e_{1,1}} + \lambda^{e_3-2e_{2,2}-e_{1,2}}
+ \lambda^{e_3-2e_{2,2}-2e_{1,1}} + \lambda^{e_3-2e_{2,2}-e_{1,1}-e_{1,2}}
+ \lambda^{e_3-2e_{2,2}-2e_{1,2}} \\
&\,\,\,\, + \lambda^{e_3-2e_{2,2}-3e_{1,1}}
+ \lambda^{e_3-2e_{2,2}-2e_{1,1}-e_{1,2}} + \lambda^{e_3-2e_{2,2}-e_{1,1}-2e_{1,2}}
+ \lambda^{e_3-2e_{2,2}-3e_{1,2}}.
\end{align*}
\end{exa}

\begin{rema}\label{rmk3.7}
{\rm We gave the formula for the character using the Demazure-type operators.
On the other hand, the character is also given by the localization formula
(\cite[Corollary 7.4]{HM}).
For example, when we set the parameters as in Example~\ref{ex3.4},
the character is computed using the localization formula as follows:}
\begin{align*}
\chi = \lambda^{e_3} &\left(\frac{1}{(1-\lambda^{-e_{1,1}})(1-\lambda^{-e_{2,1}})
(1-\lambda^{-e_{2,2}})} + \frac{\lambda^{-2e_{2,2}}}{(1-\lambda^{-e_{1,1}})
(1-\lambda^{-e_{2,1}+e_{2,2}})(1-\lambda^{e_{2,2}})} \right. \\
&\left. + \frac{\lambda^{-2e_{2,1}}}{(1-\lambda^{-e_{1,1}})
(1-\lambda^{e_{2,1}-e_{2,2}})(1-\lambda^{e_{2,1}})} +
\frac{\lambda^{-e_{1,1}}}{(1-\lambda^{e_{1,1}})(1-\lambda^{2e_{1,1}-e_{2,1}})
(1-\lambda^{-e_{1,1}-e_{2,2}})} \right. \\
&\left. + \frac{\lambda^{-3e_{1,1}-2e_{2,2}}}{(1-\lambda^{e_{1,1}})
(1-\lambda^{3e_{1,1}-e_{2,1}+e_{2,2}})(1-\lambda^{e_{1,1}+e_{2,2}})} \right. \\
&\left. + \frac{\lambda^{3e_{1,1}-2e_{2,1}}}{(1-\lambda^{e_{1,1}})
(1-\lambda^{-3e_{1,1}+e_{2,1}-e_{2,2}})(1-\lambda^{-2e_{1,1}+e_{2,1}})} \right).
\end{align*}
\end{rema}
{\rm We can check that this result agrees with the result in Example~\ref{ex3.4}.}

\vspace{0.5cm}
{\footnotesize Yuki Sugiyama \\
Department of Mathematics, Graduate School of Science and Engineering, Chuo University, \\
Kasuga, Bunkyo-Ku, Tokyo, 112--8551 Japan. \\
{\it e-mail} : y-sugi@gug.math.chuo-u.ac.jp}

\end{document}